\documentclass[11pt,a4paper,reqno]{amsart}
\usepackage{graphicx}
\usepackage{epsfig}
\usepackage{cite}
\usepackage{xcolor}
 \newtheorem{thm}{Theorem}[subsection]
 \newtheorem{cor}[thm]{Corollary}

  \newtheorem{assum}[thm]{Assumption}
  \newtheorem{pf}[thm]{Proof}
 \theoremstyle{definition}
 \newtheorem{defn}[thm]{Definition}
 \theoremstyle{remark}

\begin{document}
\setcounter{page}{1}
\begin{flushleft}
{\scriptsize }
\end{flushleft}
\bigskip
\bigskip
\title[Subjective Equilibria under Beliefs of Exogenous Uncertainty]{Subjective Equilibria under Beliefs of Exogenous Uncertainty: Linear Quadratic Case}
\author[]{G\"{u}rdal Arslan$^1$
 }
\thanks{$^1$Department of Electrical and Computer Engineering, University of Hawaii at Manoa, 2540 Dole Street, Honolulu, HI 96822, e-mail: gurdal@hawaii.edu}

\author[]{Serdar Y\"{u}ksel$^2$
 }
\thanks{$^2$Department of Mathematics and Statistics, Queen’s University, Kingston, ON K7L 3N6,
Canada, e-mail: yuksel@mast.queensu.ca}

\begin{abstract}
%Please provide an abstract of 50 to 150 words. The abstract should
%not contain any undefined abbreviations or unspecified references.

We consider a stochastic dynamic game where players have their own linear state dynamics and quadratic cost functions. Players are coupled through some environment variables, generated by another linear system driven by the states and decisions of all players. Each player observes his own states realized up to the current time as well as the past realizations of his own decisions and the environment variables.  Each player (incorrectly) believes that the environment variables are generated by an independent exogenous stochastic process. In this setup, we study the notion of ``subjective equilibrium under beliefs of exogenous uncertainty (SEBEU)'' introduced in our recent work \cite{arslan2023subjective}. At an SEBEU, each player's strategy is optimal with respect to his subjective belief; moreover, the objective probability distribution of the environment variables is consistent with players' subjective beliefs. We construct an SEBEU in pure strategies, where each player strategy is an affine function of his own state and his estimate of the system state.



\end{abstract}
\date{\today}
\maketitle

\smallskip

\section{Introduction}

In this paper, we consider a linear-quadratic non-cooperative dynamic game played in discrete-time. In our model, players have their own state and decision variables; furthermore, our model includes the so-called ``environment variables'', that are common to all players. The state and decision variables of all players as well as the environment variables belong to finite-dimensional Euclidian spaces. Each player's state  evolves according to a linear stochastic recursion driven by his own decisions as well as the environment variables . Environment variables, on the other hand, are generated by a separate linear stochastic recursion driven by the state and decision variables of all players. At each stage, players select their decisions simultaneously based on their own information available at the time, which consist of their own state variables realized up to the current time as well as the past realizations of their own decisions and the common environment variables. Following the collective decisions of all players during each stage, the environment variables for that stage are realized, players incur their stage-costs, and all of the system states evolve to the next stage. Each player's stage-cost is determined according to a quadratic function of the player's own state and decision variable as well as the common environment variables realized during the stage.
Each player wishes to minimize the expected value of his own long-term cost, that is a discounted sum of his stage-costs.

In our setup, each player's problem depends on the decisions of the other players only through the environment variables. However, a player has no knowledge about the other players, not even the presence of the other players. As mentioned above, each player observes the past realizations of the environment variables during any stage but he has no knowledge of how the environment variables are generated.  When the strategies of the other players are fixed, each player faces a partially observed sequential decision problem, however, the player is not aware of this fact. As such, each player incorrectly assumes that the environment variables, that are sequentially revealed to the player, is an independent random sequence that is exogenous to his problem.  Each player has a belief represented by a subjective probability distribution for the environment variables.
Under such beliefs, players seek to minimize their expected long-term costs by choosing ``optimal strategies''.

An equilibrium situation arises when players employ optimal strategies, with respect to their beliefs, and the resulting objective marginal probability distribution for the environment variables is consistent with their subjective beliefs.  Such strategies and beliefs are said to form a \textit{subjective equilibrium under beliefs of exogenous uncertainty (SEBEU)} \cite{arslan2023subjective}.
At an SEBEU, if the objective marginal probability distribution of the environment variables is revealed to a player (who does not question his assumption of independence and exogeneity for the environment variables)
would have no reason to change his belief; consequently, he would have no incentive to change his strategy.

The notion of SEBEU is related to the notion of price-taking behavior in mathematical economics \cite{arrow1954existence,radner1968competitive} and can be regarded as a generalization of price-taking behavior to stochastic dynamic games. The notion of SEBEU is distinct from the well-known notions of Nash equilibrium and mean-field equilibrium \cite{huang2006,lasry2007mean} as well as from the other notions of subjective equilibrium in the literature \cite{adlakha2015equilibria,weintraub2011industry,wiszniewska2017redefinition,dudebout2012empirical,dudebout2014exogenous,Kalai:GEB:1995}.
In particular, a Nash equilibrium in our setup would require each player to take into account the effects of his own decisions on the environment variables that, in turn, affects the player's cost as well as his state evolution.
In general, this leads to an infinite regression where each player needs to estimate the other players' states, the other players' estimates of their opponents' states, and so on. Players at an SEBEU ignore such closed-loop effects, which essentially decouples player problems and leads to a significant simplification of the overall problem. 
The notion of stationary equilibrium in \cite{adlakha2015equilibria} is hinged on a condition of consistency between the player beliefs on the population states being constant 
and the long-run  average of the actual population states (over an infinite-population); see also the analogous notion of oblivious equilibrium in \cite{weintraub2011industry}. We also selectively cite \cite{gmytrasiewicz2005framework} where players make optimal decisions by maintaining subjective nested beliefs over others’ types and their beliefs about others.
We refer the reader to \cite{arslan2023subjective} for a review of the related literature vis-a-vis the notion of SEBEU as well as for a motivational application on the Autonomous Demand Response (ADR) problem on price-sensitive control of energy usage in multiple (possibly very large number of) units, e.g., \cite{constantopoulos1991estia,roozbehani2012volatility}.

The existence of SEBEU and its various properties are studied in \cite{arslan2023subjective} for compact-metric space models, i.e., all state, decision, and environment variables belong to compact metric spaces. In particular, the existence of an SEBEU in \textit{mixed strategies} is established; moreover, an SEBEU is shown to be an approximate Nash equilibrium in large-scale games. In this paper, our main objective is to obtain the counterparts of the results of \cite{arslan2023subjective} for linear-quadratic models described above, which fall outside the compact-metric models considered in \cite{arslan2023subjective}. Moreover, we aim to obtain more explicit and refined results by exploiting the specific linear-quadratic structure of the models.

A significant body of literature, pioneered in large part by Tamer Ba\c{s}ar overall several seminal studies, studies multi-stage linear-quadratic games under various asymmetric information patterns and demonstrates the richness and fragility of equilibrium properties to informational changes. A linear-quadratic-Gaussian (LQG) game in which each player makes his own noisy linear observations of the system state in addition to having access to the past (but not the current) observations of the system state by his opponents (called one-step delay observation sharing pattern) is studied in \cite{ba1978two,basar1978decentralized}. The existence of a unique Nash equilibrium, that is affine in the information variables, via contraction properties of best-response dynamics is established in several studies including \cite{ba1978two,basar1978decentralized}, whereas exact expressions are presented in \cite{ba1978two} for the two-player case; see also \cite{papavassilopoulos1982linear}. A remarkable insight is that if one also shares actions, leading to one-step-delay sharing, the Nash equilibria become non-unique even when the one-step delay observation sharing would lead to a unique equilibrium \cite{ba1978two,basar1978decentralized}. For further results on informational fragility and existence of equilibria for static and dynamic problems, we refer the reader to \cite{tBasarStochasticDiffGames,basols99,sanjari2021optimality} and \cite[Chapter 10]{YukselBasarBook23}.

Among further relevant studies, \cite{gupta2014common} shows the existence of a unique Markov perfect equilibrium in two-player LQG games where the posterior beliefs of the players on the systtem state given the common information evolves independently of the player strategies. \cite{vasal2021signaling} considers two-player finite-horizon LQG games where each player observes his own state (but not his opponent's state) and the joint decisions. Under certain conditions,  \cite{vasal2021signaling} constructs a perfect Bayesian equilibrium in such games. \cite{colombino2017mutually} shows that equilibrium strategies can be computed in zero-sum two-team linear-quadratic games with decentralized information structures under a mutual quadratic invariance condition. To the best of the authors' knowledge, the existing results in the literature are not applicable to the general linear-quadratic models considered in this paper.

{\bf Contributions.}
We present sufficient conditions, which involve only the primitive problem data, under which SEBEU strategies can be explicitly constructed in \textit{pure strategies} for the linear-quadratic models. Our construction, for both finite-horizon and infinite-horizon cases, lead to pure strategies that are affine functions of players' own states and their estimates of the system state. When the primitive random variables are Gaussian, players' estimates of the system state can be generated by a linear recursion, i.e., a Kalman filter. We note that a SEBEU for the linear setup considered in our paper (when the primitive random variables are Gaussian) admits solutions that are always linear, whenever an equilibrium exists. On the other hand, a Nash equilibrium may or may not be linear and to our knowledge beyond the setups studied in \cite{ba1978two,basar1978decentralized} and noted in \cite[Chapter 10]{YukselBasarBook23}, the existence and structure of Nash equilibria in partially observed settings has not been studied in the literature. Finally, when the number of decision makers is infinite, we show that an SEBEU is also a Nash equilibrium.

{\bf Notation.} $\mathbb{N}_0$ and $\mathbb{N}$ denote the nonnegative and positive integers, respectively;  $\mathbb{R}$ denotes the real numbers; $[\cdot]^{\prime}$ denotes the transpose of a vector or matrix; $|\cdot|$ denotes the Euclidean norm for a vector; $|x|_Q^2:=x^{\prime}Qx$ for a vector $x$ and a matrix $Q$;
$A\succeq0$ denotes a nonnegative definite matrix $A$, whereas $A\succ0$ denotes a positive definite matrix $A$; 
%$A\succeq B$ denotes $A-B\succeq0$;
$\mathcal{P}(\mathbb{T})$ denotes the set of probability measures on the Borel sigma algebra $\mathcal{B}(\mathbb{T})$ of a topological space $\mathbb{T}$; $P[\cdot]$, 
$E[\cdot]$ and $\textrm{cov}[\cdot]$ denote the expectation and the covariance, respectively, ($E^{\mu}[\cdot]$ is also used to emphasize the underlying strategy $\mu$);
%$X\sim\mathcal{N}(\theta,\Sigma)$ denotes a Gaussian random vector with mean $\theta$ and covariance $\Sigma$; 
$f(x)=\mathcal{O}(x)$ means $\limsup_{x\rightarrow0} |f(x)/x|$ is finite.

\section{Subjective equilibrium under beliefs of exogenous uncertainty}
\subsection{Model and definitions}\label{se:model}
Consider a decentralized stochastic system with $N$ decision makers, $N\in\mathbb{N}$, where the $i$-th decision maker is referred to as DM$^i$, $i\in [1,N]:=\{1,\dots,N\}$. Time variable $t\in[0,T):=\{0,\dots,T-1\}$ is an integer where $T\in\mathbb{N}$ is the time horizon; we will also allow for $T= \infty$ in which case $[0,T) = \mathbb{N}_0$. Each DM$^i$ has its own state $x_t^i\in\mathbb{R}^{n_x^i}$, control input $u_t^i\in\mathbb{R}^{n_u^i}$, and random disturbance $w_t^i\in\mathbb{R}^{n_w^i}$ at time $t\in[0,T)$, where $n_x^i,n_u^i,n_w^i\in\mathbb{N}$. In addition, there are so-called environment variables denoted by $y_t\in\mathbb{R}^{n_y}$ at time $t\in[0,T)$, where $n_y\in\mathbb{N}$. Each DM$^i$ has linear state dynamics
\begin{equation}
\label{eq:x}
x_{t+1}^i  = A_t^i x_t^i +B_t^i u_t^i + C_t^i y_t + w_t^i, \qquad t\in[0,T)
\end{equation}
starting from some (possibly random) initial state $x_0^i$, and a quadratic cost function
\begin{equation}
\label{eq:c}
c_t^i(x^i,u^i,y) =|x^i|_{Q_t^i}^2 + |u^i|_{R_t^i}^2 + 2y^{\prime} (K_t^i u^i + L_t^i x^i)
\end{equation}
where $A_t^i$, $B_t^i$, $C_t^i$, $Q_t^i$, $R_t^i$, $K_t^i$, $L_t^i$ are appropriate dimensional matrices ($Q_t^i\succ0$, $R_t^i\succ0$), for all $i\in[1,N]$, $t\in[0,T)$.
Each DM$^i$ has also a quadratic terminal cost function $c_T^i(x^i) = |x^i|_{Q_T^i}^2$, where $Q_T^i\succ0$, when $T\in\mathbb{N}$ ($c_T^i\equiv0$ when $T=\infty$).
The environment variables are generated by the linear dynamics
\begin{align}
x_{t+1}^0  = & A_t^0 x_t^0 +\sum_{i\in[1,N]} (B_t^{1,i} u_t^i + B_t^{2,i} x_t^i)  + w_t^0 \label{eq:x0}\\
y_t = &  D_t x_t^0 + \sum_{i\in[1,N]} (E_t^{1,i} u_t^i + E_t^{2,i} x_t^i) +\xi_t \label{eq:y}
\end{align}
where $x_t^0\in\mathbb{R}^{n_x^0}$ is the state of the environment dynamics ($x_0^0$ is some possibly random initial state), $\xi_t\in\mathbb{R}^{n_{\xi}}$, $w_t^0\in\mathbb{R}^{n_w^0}$ are random disturbances, for some $n_x^0,n_{\xi},n_w^0\in\mathbb{N}$, and $A_t^0$, $B_t^{1,i}$, $B_t^{2,i}$, $D_t$, $E_t^{1,i}$,$E_t^{2,i}$ are appropriate dimensional matrices, for all $i\in[1,N]$, $t\in[0,T)$.

Each DM$^i$ has access to the information $(I_t^i,Y_{t-1})$ at time $t\in[0,T)$, where
\begin{equation}
%I_0^i:=x_0^i, \quad Y_{-1}:=y_{-1}, \quad 
I_t^i   := ((x_k^i)_{k\in[0,t]},(u_k^i)_{k\in[0,t)}), \quad Y_{t-1} :=(y_k)_{k\in[0,t)}, \quad t\in[0,T).
\label{eq:info}
\end{equation}
A \textit{pure strategy} for DM$^i$ is a collection $s^i=(s_t^i)_{t\in[0,T)}$  of measurable mappings $s_t^i:\mathbb{I}_t^i \times \mathbb{Y}^{[0,t)} \to \mathbb{U}^i$, where $\mathbb{I}_t^i:=(\mathbb{R}^{n_x^i})^{[0,t]}\times(\mathbb{R}^{n_u^i})^{[0,t)}$,  $\mathbb{Y}:=\mathbb{R}^{n_y}$,
$t\in[0,T)$. For both $T\in\mathbb{N}$ and $T=\infty$, $(I_T^i,Y_{T-1})\in\mathbb{I}_T^i\times\mathbb{Y}^{[0,T)}$ denotes DM$^i$'s entire history of observations.
We let $\mathbb{S}^i$ denote the set of pure strategies for each DM$^i$. If DM$^i$ employs a pure strategy $s^i\in\mathbb{S}^i$, then
\begin{equation}
\label{eq:u}
u_t^i=s_t^i(I_t^i,Y_{t-1}), \qquad t\in[0,T).
\end{equation}
For any \textit{joint pure strategy} $s=(s^1,\dots,s^N)\in\mathbb{S}:=\mathbb{S}^1\times\cdots\times\mathbb{S}^N$ employed,
each DM$^i$ incurs a long-term cost
\begin{equation}
\label{eq:ltc1}
\bar{J}^i(s^i,s^{-i}) = E \Bigg[\sum_{t\in[0,T)} (\beta^i)^t c_t^i(x_t^i,u_t^i,y_t) + (\beta^i)^T c_T^i(x_T^i) \Bigg]
\end{equation}
where $s^{-i}\in\mathbb{S}^{-i}:=\times_{j\not=i} \mathbb{S}^j$ denotes the pure strategies of all DMs other than DM$^i$,
$\beta^i\in[0,1]$ is DM$^i$'s discount factor ($\beta^i\in[0,1)$ when $T=\infty$).
The expectation in (\ref{eq:ltc1}) is taken with respect to the probability distribution over the set $\mathbb{I}_T^i\times\mathbb{Y}^{[0,T)}$ of histories  induced by $s\in\mathbb{S}$.
Note that each DM$^i$'s cost is influenced by the other DMs only through the environment variables $Y_{T-1}$.

Suppose now that each DM$^i$ views the sequence of environment variables as an {\bf independent exogenous} (that is, primitive, generated by nature independently of the all other primitive random variables) random sequence $Z_{T-1}:=(z_t)_{t\in[0,T)}$ with a given probability distribution $\zeta\in\mathcal{P}(\mathbb{Y}^{[0,T)})$  instead of the endogenous sequence $Y_{T-1}$ generated by (\ref{eq:x})-(\ref{eq:u}). We note that $Z_{T-1}$ need not be a collection of independent random variables. We define DM$^i$'s cost corresponding to a pure strategy $s^i\in\mathbb{S}^i$ and a deterministic sequence of environment variables  $\bar{Z}_{T-1}=(\bar{z}_t)_{t\in[0,T)}\in \mathbb{Y}^{[0,T)}$ as
\begin{equation}
\label{eq:ltc2}
J^i(s^i,\bar{Z}_{T-1}) := E \Bigg[\sum_{t\in[0,T)} (\beta^i)^t c_t^i\big(x_t^i,u_t^i,\bar{z}_t\big)   + (\beta^i)^T c_T^i(x_T^i) \Bigg]
\end{equation}
where
\begin{align*}
x_{t+1}^i  = & A_t^i x_t^i +B_t^i u_t^i + C_t^i \bar{z}_t + w_t^i, \qquad t\in[0,T) \\
u_t^i = & s_t^i(I_t^i,\bar{Z}_{t-1}), \qquad t\in[0,T)
\end{align*}
with $\bar{Z}_{t-1} :=(\bar{z}_k)_{k\in[0,t)}$, for $t\in[0,T)$. Note that the expectation in (\ref{eq:ltc2}) is taken with respect to the probability distribution over the set $\mathbb{I}_T^i$ of DM$^i$'s state-control histories  induced by $(s^i,\bar{Z}_{T-1})\in\mathbb{S}^i \times \mathbb{Y}^{[0,T)}$.
Accordingly, each DM$^i$ who models the environment variables as an independent exogenous random sequence with the probability distribution $\zeta\in\mathcal{P}(\mathbb{Y}^{[0,T)})$ aims to minimize $E^{\zeta} [J^i(s^i,\cdot)]$ by choosing a pure strategy $s^i\in\mathbb{S}^i$.
Therefore, from DM$^i$'s viewpoint, no DM has any influence on the environment variables and the cost $E^{\zeta} [J^i(s^i,\cdot)]$ is independent of the strategies of all DMs other than DM$^i$.

To ensure that the long-term cost (\ref{eq:ltc1}) and (\ref{eq:ltc2}) are well-defined for any joint strategy $s\in\mathbb{S}$ and independent exogenous distribution $\zeta\in\mathcal{P}(\mathbb{Y}^{[0,T)})$ for the environment variables, we make the following assumption throughout this paper.

\begin{assum}\label{as:lqudc}
$X_0, W_0,W_1,\dots,\xi_0,\xi_1,\dots,$ are mutually independent and have finite second-order moments where
$X_t:=(x_t^i)_{i\in[0,N]}$, $W_t:=(w_t^i)_{i\in[0,N]}$.
\end{assum}

\begin{defn}
\label{def:sebeu}
A joint strategy $s=(s^1,\dots,s^N)\in\mathbb{S}$ is called a {\it subjective equilibrium under beliefs of exogenous uncertainty (SEBEU) in pure strategies} if there exists subjective beliefs $\zeta^1,\dots,\zeta^N\in\mathcal{P}(\mathbb{Y}^{[0,T)})$, where $\zeta^i$ is DM$^i$'s belief about the environment variables, such that the following two conditions hold.
\begin{itemize}
\item[(i)] Each DM$^i$'s strategy $s^i$ is optimal with respect to its subjective belief $\zeta^i$, i.e.,
$$E^{\zeta^i} [J^i(s^i,\cdot)]  \leq E^{\zeta^i} [ J^i(\tilde{s}^i,\cdot) ], \quad \forall  i\in[1,N], \tilde{s}^i\in\mathbb{S}^i$$
\item[(ii)] Each DM$^i$'s subjective belief $\zeta^i$ is consistent with the objective probability distribution $\zeta_{s}$ of the environment variables, i.e.,
$$\zeta^i(B) = \zeta_s(B), \quad \forall  i\in[1,N], B\in\mathcal{B}(\mathbb{Y}^{[0,T)})$$
where $\zeta_s\in\mathcal{P}(\mathbb{Y}^{[0,T)})$ is generated by $s$ endogenously through (\ref{eq:x})-(\ref{eq:u}).\\
\end{itemize}
\end{defn}

Note that, at an SEBEU, each DM$^i$ ignores the influence of its own strategy $s^i$ on the environment variables in minimizing its long-term cost. In contrast, the well-known concept of {\it Nash equilibrium in pure strategies} requires a joint pure strategy $s=(s^1,\dots,s^N)\in\mathbb{S}$ to satisfy
\begin{equation}
\bar{J}^i(s^i,s^{-i}) \leq \bar{J}^i(\tilde{s}^i,s^{-i}), \quad \forall i\in[1,N], \tilde{s}^i\in\mathbb{S}^i
\label{eq:neq}
\end{equation}
where each DM$^i$ takes into account the entire influence of its own strategy $s^i$ on its long-term cost including through the environment variables.

%%%%%%%%%%%%%%%%%%%%%%%%%%%%%%%%%%%%%%%%%%%%%%%%%%%%%%%%%%%%%%%%%%%%%%%%%%%%%%%%
\section{Linear quadratic systems}
\label{se:lq}

In this section, we show how to construct an SEBEU in pure strategies under certain conditions.

\subsection{Finite-horizon case ($T\in\mathbb{N}$)}\label{finiteLQG}
Suppose that DM$^i$ instead aims to minimize
$E^{\zeta} [J^i(s^i,\cdot)]$ over the set of pure strategies $\mathbb{S}^i$ for an independent exogenous sequence $(z_t)_{t\in[0,T)}$ of environment variables with a given probability distribution $\zeta\in\mathcal{P}(\mathbb{Y}^{[0,T)})$. If Assumption~\ref{as:lqudc} holds and $\max_{t\in[0,T)}E[|z_t|^2]<\infty$,
an optimal strategy minimizing $E^{\zeta} [J^i(s^i,\cdot)]$ over $s^i\in\mathbb{S}^i$ is obtained as
\begin{align*}
s_t^i(I_t^i,Z_{t-1}) = & F_t^i x_t^i +  \sum_{k\in[t,T)} G_{t,k}^i E[z_k|Z_{t-1}] + H_t^i, \quad t\in[0,T)
\end{align*}
where $(F_t^i)_{t\in[0,T)}$, $(G_{t,k}^i)_{0\leq t \leq k < T}$, $(H_t^i)_{t\in[0,T)}$  are matrices of appropriate dimensions, which are pre-computable independently of the  sequence $(z_t)_{t\in[0,T)}$ of environment variables; see Appendix~\ref{appx:lqfinite}.

Therefore, a joint strategy  $s=(s^i)_{i\in[1,N]}\in\mathbb{S}$ is an SEBEU if
\begin{align}
s_t^i(I_t^i,Y_{t-1}) =
F_t^i x_t^i +  \sum_{k\in[t,T)} G_{t,k}^i E[y_k|Y_{t-1}] + H_t^i, \quad i\in[1,N], \ t\in[0,T) \label{eq:eqs}
\end{align}
and
\begin{align}
y_t  = & D_t x_t^0 + \sum_{j\in[1,N]} E_t^{1,j} \bigg(F_t^j x_t^j +  \sum_{k\in[t,T)} G_{t,k}^j E[y_k|Y_{t-1}] + H_t^j \bigg) + \sum_{j\in[1,N]} E_t^{2,j} x_t^j +\xi_t \label{eq:eqp}\\
x_{t+1}^0  = & A_t^0 x_t^0 +\sum_{j\in[1,N]} B_t^{1,j} \bigg(F_t^j x_t^j +  \sum_{k\in[t,T)} G_{t,k}^j E[y_k|Y_{t-1}] + H_t^j\bigg) +  \sum_{j\in[1,N]} B_t^{2,j} x_t^j + w_t^0 \label{eq:eqy} \\
x_{t+1}^j  = & A_t^j x_t^j +B_t^j \bigg(F_t^j x_t^j +  \sum_{k\in[t,T)} G_{t,k}^j E[y_k|Y_{t-1}] + H_t^j\bigg) + C_t^j y_t + w_t^j, \  j\in[1,N] \label{eq:eqx}
\end{align}
provided $\max_{t\in[0,T)}E[|y_t|^2]<\infty$.
This reveals that an SEBEU of the form above can be obtained if a random sequence of environment variables having finite second-order moments and satisfying (\ref{eq:eqp})-(\ref{eq:eqx}), called \textit{equilibrium environment variables}, can be constructed.
We next attempt to construct equilibrium environment variables via the conditional expectations $E[y_k|Y_{t-1}]$, $0\leq t \leq k<T$.

We write (\ref{eq:eqp})-(\ref{eq:eqx}) in the vector form
\begin{align}
y_t  = & \mathcal{D}_t X_t + \sum_{k\in[t,T)} \mathcal{G}_{t,k}^p E[y_k|Y_{t-1}] + \mathcal{H}_t^p  +\xi_t \label{eq:eqpv1}\\
X_{t+1}  = & \mathcal{A}_t X_t + \sum_{k\in[t,T)} \mathcal{G}_{t,k}^X E[y_k|Y_{t-1}] + \mathcal{H}_t^X + \mathcal{C}_t y_t + W_t \label{eq:eqxv1}
\end{align}
where $\mathcal{D}_t$, $\mathcal{G}_{t,k}^p$, $\mathcal{H}_t^p$, $\mathcal{A}_t$, $\mathcal{G}_{t,k}^X$, $\mathcal{H}_t^X$, $\mathcal{C}_t$, $0\leq t \leq k < T$, are matrices of appropriate dimension.
For each $k\in[0,T)$ and arbitrary $\hat{X}_{k|k-1}$, we introduce the linear equations
\begin{equation}
\Bigg\{
\begin{array}{rl}
\hat{y}_{t|k-1}  = & \mathcal{D}_t \hat{X}_{t|k-1} + \sum_{n\in[t,T)} \mathcal{G}_{t,n}^p \hat{y}_{n|k-1} + \mathcal{H}_t^p  + \hat{\xi}_t \\
\hat{X}_{t+1|k-1}  = & \mathcal{A}_t \hat{X}_{t|k-1} + \sum_{n\in[t,T)} \mathcal{G}_{t,n}^X \hat{y}_{n|k-1} + \mathcal{H}_t^X + \mathcal{C}_t \hat{y}_{t|k-1} + \hat{W}_t
\end{array} \Bigg\}, \ t\in[k,T) \label{eq:eqpxc1}
\end{equation}
where $\hat{\xi}_t:=E[\xi_t]$, $\hat{W}_t:=E[W_t]$, and $\hat{y}_{t|k-1}$ are the unknowns representing the conditional expectations $E[y_t|Y_{k-1}]$, $t\in[k,T)$.
We will refer to (\ref{eq:eqpxc1}) as $Eq(k,\hat{X}_{k|k-1})$.
Evidently, the existence of equilibrium environment variables $(y_t)_{t\in[0,T)}$ requires solutions $\hat{y}_{t|k-1}$ to $Eq(k,E[X_k|Y_{k-1}])$ for all realizations of $Y_{k-1}$ that can be generated by (\ref{eq:eqpv1})-(\ref{eq:eqxv1}).\\

\begin{assum}
\label{as:1}
For each $k\in[0,T)$ and arbitrary $\hat{X}_{k|k-1}$, the linear equations $Eq(k,\hat{X}_{k|k-1})$ in (\ref{eq:eqpxc1}) has a \textit{unique} solution. We note that such a  solution is necessarily in the form
\begin{equation}
	\hat{y}_{t|k-1} = a_{t,k-1} \hat{X}_{k|k-1} + b_{t,k-1}, \qquad t\in[k,T) \nonumber
\end{equation}
where $(a_{t,k-1})_{t\in[k,T)}$ and $(b_{t,k-1})_{t\in[k,T)}$ depend only on the constants $\mathcal{D}_t$, $\mathcal{G}_{t,n}^{p}$, $\mathcal{H}_t^p$, $\hat{\xi}_t$, $\mathcal{A}_t$, $\mathcal{G}_{t,n}^X$, $\mathcal{H}_t^X$, $\mathcal{C}_t$, $\hat{W}_t$, $0\leq k \leq t \leq n < T$.
\end{assum}

Assumption~\ref{as:1} is satisfied, for example, when the matrices $E_t^{1,j}$, $E_t^{2,j}$, $B_t^{1,j}$, $B_t^{2,j}$, for all $j\in[1,N]$, $t\in[0,T)$, are sufficiently small in size, i.e., when the sensitivity of the environment variables to the states and decisions of DMs is sufficiently small. Note that the linear equations $Eq(k,\hat{X}_{k|k-1})$, for any $k\in[0,T)$, can be written as $$\left[\begin{array}{c}\hat{y}_{k|k-1} \\ \vdots \\ \hat{y}_{T-1|k-1} \end{array}\right] = \Lambda_k \left[\begin{array}{c}\hat{y}_{k|k-1} \\ \vdots \\ \hat{y}_{T-1|k-1} \end{array}\right]  + \Upsilon_k \hat{X}_{k|k-1}$$ where $\Lambda_k$, $\Upsilon_k$ are constants of appropriate dimension, which depend on $\mathcal{D}_t$, $\mathcal{G}_{t,n}^{p}$, $\mathcal{H}_t^p$, $\hat{\xi}_t$, $\mathcal{A}_t$, $\mathcal{G}_{t,n}^X$, $\mathcal{H}_t^X$, $\mathcal{C}_t$, $\hat{W}_t$, $0\leq k \leq t \leq n < T$. Furthermore,
the matrices $F_t^j$, $G_{t,k}^j$, $H_t^j$ defining each DM$^j$'s optimal controller are independent of $E_t^{1,j}$, $E_t^{2,j}$, $B_t^{1,j}$, $B_t^{2,j}$. Accordingly, when the matrices $E_t^{1,j}$, $E_t^{2,j}$, $B_t^{1,j}$, $B_t^{2,j}$, for all $j\in[1,N]$, $t\in[0,T)$, are sufficiently small in size, $I-\Lambda_k$ is invertible.

Under Assumption~\ref{as:lqudc} and \ref{as:1}, it is natural to introduce the following \textit{candidate} equilibrium environment variables, which is
an independent exogenous sequence $(\bar{y}_t)_{t\in[0,T)}$ generated by the recursion,
for $t\in[0,T)$,
\begin{align}
\bar{y}_t  = & \mathcal{D}_t \bar{X}_t + \sum_{n\in[t,T)} \mathcal{G}_{t,n}^{p} (a_{n,t-1} E[\bar{X}_t|\bar{Y}_{t-1}] + b_{n,t-1}) + \mathcal{H}_t^p +\xi_t \label{eq:eqpbar1}\\
\bar{X}_{t+1}  = & \mathcal{A}_t \bar{X}_t + \sum_{n\in[t,T)} \mathcal{G}_{t,n}^{X} (a_{n,t-1} E[\bar{X}_t|\bar{Y}_{t-1}] + b_{n,t-1}) + \mathcal{H}_t^X+ \mathcal{C}_t \bar{y}_t  + W_t \label{eq:eqXbar1}
\end{align}
where $(a_{n,t-1})_{n\in[t,T)}$ and $(b_{n,t-1})_{n\in[t,T)}$ are as in Assumption~\ref{as:1}, $\bar{Y}_{t-1} = (\bar{y}_k)_{k\in[0,t)}$, for $t\in[0,T)$, and
 $\bar{X}_0=X_0$.
 It can be shown that the random sequence $(\bar{y}_t)_{t\in[0,T)}$ generated by (\ref{eq:eqpbar1})-(\ref{eq:eqXbar1}) is indeed a sequence of equilibrium environment variables. This leads to the existence of SEBEU.\\

\begin{thm}
\label{thm:ex}
Consider the finite-horizon linear quadratic problem of this subsection.
If Assumption~\ref{as:lqudc} and \ref{as:1} hold,
there exists an SEBEU strategy $s\in\mathbb{S}$ defined by
\begin{align}
s_t^i(I_t^i,Y_{t-1}) =F_t^i x_t^i +G_{t}^i \hat{\bar{X}}_{t|t-1}(Y_{t-1}) + \check{H}_t^i, \qquad  i\in[1,N], \ t\in[0,T)
\label{eq:eqstrf}
\end{align}
where $G_t^i := \sum_{k\in[t,T)} G_{t,k}^i a_{k,t-1}$, $\check{H}_t^i := \sum_{k\in[t,T)} G_{t,k}^i b_{k,t-1} + H_t^i$, and
$\hat{\bar{X}}_{t|t-1}(\cdot)$ denotes the state estimator mapping from $\bar{Y}_{t-1}$ to $E[\bar{X}_t|\bar{Y}_{t-1}]$ based on (\ref{eq:eqpbar1})-(\ref{eq:eqXbar1}).
%If there are sets of identical DMs, there exists a symmetric SEBEU $s\in\mathbb{S}$ such that identical DMs have identical strategies at  $s$, i.e., $s^i=s^j$, if DM$^i$ and DM$^j$ are identical.
\\
\end{thm}

\begin{pf}
Consider any $0\leq k\leq\ell< T$ and arbitrary $\hat{X}_{k|k-1}$.
The equations $Eq(k,\hat{X}_{k|k-1})$  for $t\in[\ell,T)$ coincide with the equations $Eq(\ell,\hat{X}_{\ell|k-1})$, where $\hat{X}_{\ell|k-1}$ is obtained by propagating $\hat{X}_{k|k-1}$ through (\ref{eq:eqpxc1}). By Assumption~\ref{as:1}, this implies that, for $t\in[\ell,T)$,
\begin{equation}
a_{t,k-1} \hat{X}_{k|k-1} + b_{t,k-1} = a_{t,\ell-1} \hat{X}_{\ell|k-1} + b_{t,\ell-1}.
\label{eq:kel}
\end{equation}

By Assumption~\ref{as:1} and substituting (\ref{eq:kel}) into $Eq(k,\hat{X}_{k|k-1})$, we obtain, for $t\in[k,T)$,
\begin{align}
\hat{y}_{t|k-1}  = & \mathcal{D}_t \hat{X}_{t|k-1} + \sum_{n\in[t,T)} \mathcal{G}_{t,n}^p (a_{n,t-1} \hat{X}_{t|k-1} + b_{n,t-1}) + \mathcal{H}_t^p  + \hat{\xi}_t \label{eq:eqp3}\\
\hat{X}_{t+1|k-1}  = & \mathcal{A}_t \hat{X}_{t|k-1} + \sum_{n\in[t,T)} \mathcal{G}_{t,n}^X (a_{n,t-1} \hat{X}_{t|k-1} + b_{n,t-1}) + \mathcal{H}_t^X + \mathcal{C}_t \hat{y}_{t|k-1} + \hat{W}_t \label{eq:eqX3}
\end{align}
where $\hat{y}_{t,k-1}=a_{t,k-1} \hat{X}_{k|k-1} + b_{t,k-1}$.
Taking conditional expectations of both sides of (\ref{eq:eqpbar1})-(\ref{eq:eqXbar1}) given $\bar{Y}_{k-1}$ leads to, for $t\in[k,T)$,
\begin{align}
E[\bar{y}_t|\bar{Y}_{k-1}]  = & \mathcal{D}_t E[\bar{X}_t|\bar{Y}_{k-1}] + \sum_{n\in[t,T)} \mathcal{G}_{t,n}^{p} (a_{n,t-1} E[\bar{X}_t|\bar{Y}_{k-1}] + b_{n,t-1}) + \mathcal{H}_t^p +\hat{\xi}_t \label{eq:eqpbar2}\\
E[\bar{X}_{t+1}|\bar{Y}_{k-1}]  = & \mathcal{A}_t E[\bar{X}_t|\bar{Y}_{k-1}] + \sum_{n\in[t,T)} \mathcal{G}_{t,n}^{X} (a_{n,t-1} E[\bar{X}_t|\bar{Y}_{k-1}] + b_{n,t-1}) + \mathcal{H}_t^X  + \mathcal{C}_t E[\bar{y}_t|\bar{Y}_{k-1}]  + \hat{W}_t. \label{eq:eqXbar2}
\end{align}
Equations (\ref{eq:eqp3})-(\ref{eq:eqX3}) hold, in particular, for $\hat{X}_{k|k-1}=E[\bar{X}_k|\bar{Y}_{k-1}]$, in which case $\hat{y}_{t,k-1}= a_{t,k-1} E[\bar{X}_k|\bar{Y}_{k-1}] + b_{t,k-1} $.
Therefore, setting $\hat{X}_{k|k-1}=E[\bar{X}_k|\bar{Y}_{k-1}]$ in (\ref{eq:eqp3})-(\ref{eq:eqX3}) and comparing the right-hand-sides of (\ref{eq:eqp3})-(\ref{eq:eqX3}) and (\ref{eq:eqpbar2})-(\ref{eq:eqXbar2}) for $t=k$  yields
\begin{align*}
E[\bar{y}_k|\bar{Y}_{k-1}]  = a_{k,k-1} E[\bar{X}_k|\bar{Y}_{k-1}] + b_{k,k-1} \qquad\textrm{and} \qquad
E[\bar{X}_{k+1}|\bar{Y}_{k-1}]  = \hat{X}_{k+1|k-1}.
\end{align*}
Using $E[\bar{X}_{k+1}|\bar{Y}_{k-1}]  = \hat{X}_{k+1|k-1}$ and comparing the right-hand-sides of (\ref{eq:eqp3})-(\ref{eq:eqX3}) and (\ref{eq:eqpbar2})-(\ref{eq:eqXbar2}) for $t=k+1$  yields
\begin{align*}
E[\bar{y}_{k+1}|\bar{Y}_{k-1}]  = a_{k+1,k-1} E[\bar{X}_{k}|\bar{Y}_{k-1}] + b_{k+1,k-1} \qquad\textrm{and} \qquad
E[\bar{X}_{k+2}|\bar{Y}_{k-1}]  = \hat{X}_{k+2|k-1}.
\end{align*}
Using $E[\bar{X}_{k+2}|\bar{Y}_{k-1}]  = \hat{X}_{k+2|k-1}$ and continuing in a similar fashion yields
$$E[\bar{y}_t|\bar{Y}_{k-1}]=a_{t,k-1} E[\bar{X}_k|\bar{Y}_{k-1}] + b_{t,k-1}, \qquad t\in[k,T).$$ Using this in (\ref{eq:eqpbar1})-(\ref{eq:eqXbar1}) shows that the environment variables $(\bar{y}_t)_{t\in[0,T)}$ generated by (\ref{eq:eqpbar1})-(\ref{eq:eqXbar1})  satisfy (\ref{eq:eqpv1})-(\ref{eq:eqxv1}).
Furthermore, due to Assumption~\ref{as:1}, $(\bar{y}_t)_{t\in[0,T)}$ have finite second-order moments. This implies that $(\bar{y}_t)_{t\in[0,T)}$ are equilibrium environment variables, and the strategy defined by (\ref{eq:eqstrf}) is an SEBEU.
%Furthermore, the unique solvability of the equations $Eq(k,\hat{X}_{k|k-1})$ in Assumption~\ref{as:1} implies that equilibrium environment variables have a unique distribution, given by distribution of the environment variables generated by (\ref{eq:eqpbar1})-(\ref{eq:eqXbar1}).
%The last part follows from the fact that identical DMs have identical optimal responses to any independent exogenous sequence of environment variables.
\end{pf}

We remark that Assumption~\ref{as:1} requires the existence of a unique solution to $Eq(k,\hat{X}_{k|k-1})$ for arbitrary $\hat{X}_{k|k-1}$ which is in general not necessary for the existence of equilibrium environment variables. Assumption~\ref{as:1} can be relaxed so that $Eq(k,\hat{X}_{k|k-1})$ has a solution only for those $\hat{X}_{k|k-1}=E[\bar{X}_k|\bar{Y}_{k-1}]$ generated by the dynamics (\ref{eq:eqpbar1})-(\ref{eq:eqXbar1}). However, by also noting that an SEBEU must have the form (\ref{eq:eqs}) except on a set of observation histories of probability zero, we conclude that Assumption~\ref{as:1} cannot be relaxed in a significant way.

At an SEBEU, DMs model the environment variables as an independent exogenous stochastic process  as in (\ref{eq:eqpbar1})-(\ref{eq:eqXbar1}) where $\bar{X}_t$ is a replica of the system state $X_t$. This model for the environment variables and the resulting distribution of the environment variables is consistent with the actual environment variables generated endogenously by the equilibrium strategies. Therefore, one can start with such a model for the environment variables (instead of an arbitrary distribution for the environment variables) and establish the existence of an SEBEU in an alternative fashion. 

In the next subsection, we extend this result to infinite-horizon linear quadratic problems.

\subsection{Infinite-horizon case ($T=\infty$)}\label{ss:lqgfp}
In this subsection, we further assume that the system is time-invariant. Hence, the matrices $A_t^i$, $B_t^i$, $C_t^i$, $A_t^0$, $B_t^{1,i}$, $B_t^{2,i}$, $D_t$, $E_t^{1,i}$, $E_t^{2,i}$, $Q_t^i$, $R_t^i$, $K_t^i$, $L_t^i$ are all assumed to be independent of the time variable $t$ (hence the subscript $t$ is dropped from these matrices). Recall that $\beta^i\in[0,1)$, for all $i\in[1,N]$, in the infinite-horizon case.

\begin{assum}
\label{as:lqudcinf}
$\mbox{}$

\begin{itemize}
\item[(i)] $A^0$ is stable, $(A^i,B^i)$ is stabilizable, for all $i\in[1,N]$
\item[(ii)]  $X_0$, $W_0,W_1,\dots$, $\xi_0,\xi_1,\dots$  are mutually independent and have finite second-order moments 
\item[(iii)] $(W_t^0)_{t\in\mathbb{N}_0}$, $(\xi_t)_{t\in\mathbb{N}_0}$ are each independent and identically distributed (iid).
\end{itemize}
\end{assum}

Consider an independent exogenous sequence $(z_t)_{t\in[0,\infty)}$ of environment variables with a given probability distribution $\zeta\in\mathcal{P}(\mathbb{Y}^{[0,\infty)})$ satisfying $\sup_{t\in[0,\infty)} E[|z_t|^2]<\infty$.
An optimal strategy minimizing $E^{\zeta} [J^i(s^i,\cdot)]$ over $s^i\in\mathbb{S}^i$ is given by
$$
u_t^i =s_t^i(I_t^i,Z_{t-1})= F^i x_t^i  +\sum_{n\in\mathbb{N}_0} G_n^i E[z_{t+n}|Z_{t-1}] + H^i, \qquad t\in\mathbb{N}_0
$$
where $F^i$, $(G_n^i)_{n\in\mathbb{N}_0}$, $H^i$ are constants of appropriate dimensions, which are pre-computable independently of the sequence $(z_t)_{t\in[0,\infty)}$ of environment variables; see Appendix~\ref{appx:lqinfinite}.

Consider now a joint strategy $s=(s^i)_{i\in[1,N]}\in\mathbb{S}$ where the sequence $(y_t)_{t\in[0,\infty)}$ of environment variables generated under $s$ has the probability distribution $\zeta_s$ and $\sup_{t\in[0,\infty)}E[|y_t|^2]<\infty$. The joint strategy $s$ is an SEBEU if, for each $i\in[1,N]$, $s^i$ minimizes $E^{\zeta_s} [J^i(s^i,\cdot)]$, e.g.,
\begin{align}
s_t^i(I_t^i,Y_{t-1}) =F^i x_t^i+\sum_{n\in\mathbb{N}_0} G_n^i E[y_{t+n}|Y_{t-1}] + H^i, \quad t\in\mathbb{N}_0  \label{eq:eqstr}
\end{align}
where
\begin{align}
y_t  = & D x_t^0 +  \sum_{j\in[1,N]} E^{1,j} \bigg(F^j x_t^j+\sum_{n\in\mathbb{N}_0} G_n^j E[y_{t+n}|Y_{t-1}] + H^j\bigg) +  \sum_{j\in[1,N]} E^{2,j}  x_t^j+\xi_t \label{eq:eqp2}\\
x_{t+1}^0  = & A^0 x_t^0 + \sum_{j\in[1,N]} B^{1,j} \bigg(F^j x_t^j+\sum_{n\in\mathbb{N}_0} G_n^j E[y_{t+n}|Y_{t-1}] + H^j\bigg)  +  \sum_{j\in[1,N]} B^{2,j}  x_t^j+ w_t^0 \label{eq:eqy2} \\
x_{t+1}^j  = & A^j x_t^j +B^j \bigg(F^j x_t^j+\sum_{n\in\mathbb{N}_0} G_n^j E[y_{t+n}|Y_{t-1}] + H^j\bigg) + C^j y_t + w_t^j, \  j\in[1,N]. \label{eq:eqx2}
\end{align}
Therefore, an SEBEU can be obtained via (\ref{eq:eqstr}) provided a sequence $(y_t)_{t\in[0,\infty)}$ of environment variables satisfying (\ref{eq:eqp2})-(\ref{eq:eqx2}) and $\sup_{t\in[0,\infty)} E[|y_t|^2] < \infty$, called \textit{equilibrium environment variables}, exists. Next, we introduce a set of assumptions to ensure the existence of equilibrium environment variables.

We rewrite (\ref{eq:eqp2})-(\ref{eq:eqx2}) in the vector form
\begin{align}
y_t  = & \mathcal{D} X_t + \sum_{n\in\mathbb{N}_0} \mathcal{G}_n^{p} E[y_{t+n}|Y_{t-1}]  +\mathcal{H}^p +\xi_t \label{eq:eqpv}\\
X_{t+1}  = & \mathcal{A} X_t + \sum_{n\in\mathbb{N}_0} \mathcal{G}_n^X E[y_{t+n}|Y_{t-1}] + \mathcal{H}^X+ \mathcal{C} y_t  + W_t \label{eq:eqxv}
\end{align}
where $\mathcal{D}$, $(\mathcal{G}_n^p)_{n\in\mathbb{N}_0}$, $\mathcal{H}^p$,  $\mathcal{A}$, $(\mathcal{G}_n^X)_{n\in\mathbb{N}_0}$, $\mathcal{H}^X$,  $\mathcal{C}$ are some constants of appropriate dimensions.
Taking the expectation of both sides of (\ref{eq:eqpv})-(\ref{eq:eqxv}) shows that the existence of equilibrium environment variables requires a solution $(\hat{y}_t)_{t\in\mathbb{N}_0}$ to the infinite system of linear equations, for $t\in\mathbb{N}_0$,
\begin{align}
\hat{y}_{t}  = & \mathcal{D} \hat{X}_{t} + \sum_{n\in\mathbb{N}_0} \mathcal{G}_n^{p} \hat{y}_{t+n}  +\mathcal{H}^p +\hat{\xi}_0 \label{eq:hatpm1}\\
\hat{X}_{t+1}  = & \mathcal{A} \hat{X}_{t} + \sum_{n\in\mathbb{N}_0} \mathcal{G}_n^X \hat{y}_{t+n} + \mathcal{H}^X+ \mathcal{C} \hat{y}_{t}   + \hat{W}_0 \label{eq:hatXm1}
\end{align}
where $\hat{X}_{0}=E[X_0]$.
Accordingly, we make the following assumption.\\

\begin{assum}
\label{as:2}
The infinite system of linear equations (\ref{eq:hatpm1})-(\ref{eq:hatXm1}) has a \textit{unique} solution $(\hat{y}_{t})_{t\in\mathbb{N}_0}$ for \textit{any} $\hat{X}_{0}$, which will necessarily be in the form
\begin{equation}
\hat{y}_{t} = a_t \hat{X}_{0} + b_t, \qquad t\in\mathbb{N}_0
\label{eq:eqcp}
\end{equation}
where $(a_t)_{t\in\mathbb{N}_0}$, $(b_t)_{t\in\mathbb{N}_0}$ are assumed to be uniformly bounded and depend only on the constants $\mathcal{D}$, $(\mathcal{G}_n^p)_{n\in\mathbb{N}_0}$, $\mathcal{H}^p$, $\hat{\xi}_0$, $\mathcal{A}$, $(\mathcal{G}_n^X)_{n\in\mathbb{N}_0}$, $\mathcal{H}^X$, $\mathcal{C}$, $\hat{W}_0$.\\
\end{assum}

For any $k\in\mathbb{N}_0$, taking the conditional expectation of both sides of (\ref{eq:eqpv})-(\ref{eq:eqxv}) given $Y_{k-1}$, where $t\in[k,\infty)$, also leads to  (\ref{eq:hatpm1})-(\ref{eq:hatXm1}) with $\hat{X}_0$ replaced with $E[X_k|Y_{k-1}]$.
Therefore, under Assumption~\ref{as:2}, the conditional expectations of equilibrium environment variables satisfying (\ref{eq:eqpv})-(\ref{eq:eqxv}), if exist, are given by
\begin{equation}
E[y_{n+k}|Y_{k-1}] = a_n E[X_k|Y_{k-1}] + b_n, \qquad n\in\mathbb{N}_0, \ k\in\mathbb{N}_0.
\label{eq:hatptm1}
\end{equation}
This prompts us to introduce the independent exogenous environment variables $(\bar{y}_t)_{t\in[0,\infty)}$ generated by the recursion, for $t\in\mathbb{N}_0$,
\begin{align}
\bar{y}_t  = & \mathcal{D} \bar{X}_t + \mathcal{G}^{p} E[\bar{X}_t|\bar{Y}_{t-1}] + \check{\mathcal{H}}^p +\bar{\xi}_t \label{eq:eqpvv}\\
\bar{X}_{t+1}  = & \mathcal{A} \bar{X}_t + \mathcal{G}^X E[\bar{X}_t|\bar{Y}_{t-1}] + \check{\mathcal{H}}^X+ \mathcal{C} \bar{y}_t  + \bar{W}_t \label{eq:eqxvv}
\end{align}
 where
 $\mathcal{G}^{p}:=\sum_{n\in\mathbb{N}_0} \mathcal{G}_n^{p} a_n$, $\check{\mathcal{H}}^p := \sum_{n\in\mathbb{N}_0} \mathcal{G}_n^{p} b_n + \mathcal{H}^p$, $\mathcal{G}^{X}:=\sum_{n\in\mathbb{N}_0} \mathcal{G}_n^{X} a_n$, $\check{\mathcal{H}}^X := \sum_{n\in\mathbb{N}_0} \mathcal{G}_n^{X} b_n + \mathcal{H}^X$,
 the distribution of 
 $(\bar{X}_0,(\bar{\xi})_{t\in\mathbb{N}_0},(\bar{W})_{t\in\mathbb{N}_0})$ is identical to that of $(X_0,(\xi)_{t\in\mathbb{N}_0},(W)_{t\in\mathbb{N}_0})$, and $\bar{Y}_{t-1}=(\bar{y}_k)_{k\in[0,t)}$. By construction, the environment variables  $(\bar{y}_t)_{t\in[0,\infty)}$ generated by (\ref{eq:eqpvv})-(\ref{eq:eqxvv}) satisfy (\ref{eq:eqpv})-(\ref{eq:eqxv}); hence the environment variables $(\bar{y}_t)_{t\in[0,\infty)}$ are equilibrium environment variables, provided $\sup_{t\in[0,\infty)}E[|\bar{y}_t|^2]<\infty$. \\

\begin{thm}
\label{thm:eqexinf}
Consider the infinite-horizon linear quadratic problem of this subsection under Assumption~\ref{as:lqudcinf}-\ref{as:2}. If $(\bar{y}_t)_{t\in[0,\infty)}$ generated by (\ref{eq:eqpvv})-(\ref{eq:eqxvv}) satisfies $\sup_{t\in[0,\infty)} E[|\bar{y}_t|^2]<\infty$,
there exists an SEBEU strategy $s\in\mathbb{S}$ defined by
\begin{align}
s_t^i(I_t^i,Y_{t-1}) =F^i x_t^i+G^i \hat{\bar{X}}_{t|t-1}(Y_{t-1}) + \check{H}^i, \qquad  i\in[1,N], \ t\in\mathbb{N}_0
\label{eq:etaeq}
\end{align}
where $G^i:=\sum_{n\in\mathbb{N}_0} G_n^i  a_n$, $\check{H}^i:=\sum_{n\in\mathbb{N}_0} G_n^i  b_n + H^i$, and
$\hat{\bar{X}}_{t|t-1}(\cdot)$ denotes the state estimator mapping $\bar{Y}_{t-1} \mapsto E[\bar{X}_t|\bar{Y}_{t-1}]$ based on (\ref{eq:eqpvv})-(\ref{eq:eqxvv}).
%If there are sets of identical DMs, there exists a symmetric SEBEU $s\in\mathbb{S}$ such that identical DMs have identical strategies at  $s$, i.e., $s^i=s^j$, if DM$^i$ and DM$^j$ are identical.
\\
\end{thm}

\begin{proof}
The assumptions and the discussion preceding the theorem ensures that the environment variables generated by the joint strategy $s$ defined by (\ref{eq:etaeq}) are equilibrium environment variables (which are distributed identically to the independent exogenous environment variables generated by (\ref{eq:eqpvv})-(\ref{eq:eqxvv})).
%
%The last part follows from the fact that identical DMs have identical optimal responses to any independent exogenous sequence of environment variables. This proves the theorem.
\end{proof}

Theorem~\ref{thm:eqexinf} constructs an SEBEU under the unique solvability of (\ref{eq:hatpm1})-(\ref{eq:hatXm1}) for any $\hat{X}_0$ and the uniform boundedness of the second moments of the environment variables. These conditions cannot be relaxed as long as we insist on the existence of an equilibrium (with environment variables having uniformly bounded second-order moments) for every initial state distribution.
The uniform boundedness of the second moments of the environment variables generated under (\ref{eq:etaeq}) can be ensured if the closed-loop dynamics (identical to (\ref{eq:eqpvv})-(\ref{eq:eqxvv})) are stable. The closed-loop stability will be discussed in more detail in the Gaussian case shortly.\\

\subsection{Gaussian case} We now turn our attention to the Gaussian case where all primitive random variables are jointly Gaussian. With this Gaussian assumption, in both finite-horizon and infinite-horizon cases, the state estimator $\hat{\bar{X}}_{t|t-1}(\cdot)$ mapping $\bar{Y}_{t-1} \mapsto E[\bar{X}_t | \bar{Y}_{t-1}]$can be implemented as the following Kalman filter, for $t\in[0,T)$,
\begin{align}
\hat{\bar{X}}_{t+1|t}  = & \mathcal{A}_t \hat{\bar{X}}_{t|t-1} + \mathcal{G}_t^X  \hat{\bar{X}}_{t|t-1} + \check{\mathcal{H}}_t^X+ \mathcal{C}_t \bar{y}_t  + \hat{W}_t \nonumber
\\ & + \mathcal{A}_t\Sigma_{t|t-1} \mathcal{D}_t^{\prime} (\mathcal{D}_t\Sigma_{t|t-1}\mathcal{D}_t^{\prime}+\mbox{cov}[\xi_t])^{-1} (\bar{y}_t-(\mathcal{D}_t + \mathcal{G}_t^{p})  \hat{\bar{X}}_{t|t-1} - \check{\mathcal{H}}_t^p - \hat{\xi}_t) \label{eq:kfX}\\
\Sigma_{t+1|t}  = &
\mathcal{A}_t\big(\Sigma_{t|t-1} - \Sigma_{t|t-1} \mathcal{D}_t^{\prime} (\mathcal{D}_t\Sigma_{t|t-1}\mathcal{D}_t^{\prime}+\mbox{cov}[\xi_t])^{-1} \mathcal{D}_t  \Sigma_{t|t-1} \big) \mathcal{A}_t^{\prime}
 + \mbox{cov}[W_t] \label{eq:kfS}
\end{align}
where $\bar{X}_{0|-1}=E[X_0]$, $\Sigma_{0|-1} := \mbox{cov}[\bar{X}_0]$ (assuming $\mathcal{D}_t\Sigma_{t|t-1}\mathcal{D}_t^{\prime}+\mbox{cov}[\xi_t]$ is invertible for $t\in[0,T)$).
Therefore, under the Gaussian assumption, players' decisions at an SEBEU can be generated by a linear time-varying system in both finite-horizon and infinite-horizon cases. 

In the case of infinite-horizon, it is possible to obtain an SEBEU at which players' decisions can be generated by a linear time-invariant system. We recall that, in this case, 
$\mathcal{A}_t$, $\mathcal{C}_t$, $\mathcal{D}_t$, $\mathcal{G}_t^X$, $\check{\mathcal{H}}_t^X$, $\mathcal{G}_t^{p}$, $\check{\mathcal{H}}_t^p$, $\hat{W}_t$,  $\mbox{cov}[W_t]$,  $\hat{\xi}_t$, $\mbox{cov}[\xi_t]$ are independent of $t$.
Hence, a time-invariant SEBEU can obtained if $\mbox{\rm cov}[X_0] = \Sigma_{0|-1}=\Sigma$ where $\Sigma$ satisfies
\begin{align}
%\label{eq:ssi1}
%\hat{\bar{X}}_0  = & (\mathcal{A} + \mathcal{G}^X + \mathcal{C} (\mathcal{D} + \mathcal{G}^p)) \hat{\bar{X}}_0 + \check{\mathcal{H}}^X + \mathcal{C} ( \check{\mathcal{H}}^p +\hat{\xi}_0)  + \hat{W}%_0 \\
\label{eq:ssi2}
\Sigma  = & \mathcal{A}\big(\Sigma - \Sigma \mathcal{D}^{\prime} (\mathcal{D}\Sigma\mathcal{D}^{\prime}+\mbox{cov}[\xi_0])^{-1} \mathcal{D}  \Sigma \big) \mathcal{A}^{\prime} + \mbox{cov}[W_0].
\end{align}
This ensures $\Sigma_{t|t-1}=\Sigma$, for all $t\in\mathbb{N}_0$ and leads to the following corollary to Theorem~\ref{thm:eqexinf}.\\

\begin{cor}
\label{cor:eqexinf}
Consider the infinite-horizon linear quadratic problem of subsection~\ref{ss:lqgfp} under the assumptions of Theorem~\ref{thm:eqexinf}. Assume further that (i) $\big(X_0,(\xi_t)_{t\in\mathbb{N}_0},(W_t)_{t\in\mathbb{N}_0}\big)$ is Gaussian, (ii) $\mathcal{D}\Sigma\mathcal{D}+\mbox{\rm cov}[\xi_0]$ is invertible, and (iii) $\mbox{\rm cov}[X_0] = \Sigma$ where $\Sigma$ solves  (\ref{eq:ssi2}).
The SEBEU strategy $s$ in Theorem~\ref{thm:eqexinf} can be implemented by the linear-time invariant system: $i\in[1,N]$, $t\in\mathbb{N}_0$,
\begin{align*}
s_t^i(I_t^i,Y_{t-1}) = & F^i x_t^i+G^i \hat{\bar{X}}_{t|t-1} + \check{H}^i  \\
\hat{\bar{X}}_{t+1|t}  = & \mathcal{A} \hat{\bar{X}}_{t|t-1} + \mathcal{G}^X  \hat{\bar{X}}_{t|t-1} + \check{\mathcal{H}}^X+ \mathcal{C} y_t  + \hat{W}_0 \nonumber
\\ & + \mathcal{A}\Sigma \mathcal{D}^{\prime} (\mathcal{D}\Sigma \mathcal{D}^{\prime}+\mbox{\rm cov}[\xi_0])^{-1} (y_t-(\mathcal{D} + \mathcal{G}^{p})  \hat{\bar{X}}_{t|t-1} - \check{\mathcal{H}}^p - \hat{\xi}_0).
\end{align*}
(provided $(y_t)_{t\in[0,\infty)}$ generated under $s$ satisfies $\sup_{t\in[0,\infty)} E[|y_t|^2]<\infty$). \\
\end{cor}

Implicit in the steady-state initialization condition  $\mbox{\rm cov}[X_0] = \Sigma$ is that there exist  $\Sigma\succeq0$  solving (\ref{eq:ssi2}), which can be established under reasonable assumptions ensuring closed-loop stability.
Note that, since $A^0$ is stable, $\mathcal{A}$ would be stable if $A^i+B^iF^i$, $i\in[1,N]$ are stable, for all $i\in[1,N]$. In other words, $\mathcal{A}$ would be stable if each DM's own state dynamics (as the environment variables are seen as external inputs) is stabilized by its own control strategy.
However, the internal dynamics of the environment variables and each DM's own state dynamics are coupled in a feedback loop. The strategies of DMs designed to optimize their individual objectives by viewing the process of environment variables as independent and exogenous cannot be expected to stabilize this interconnected feedback system. The overall stability of such an interconnected system (of stable subsystems) can be guaranteed under a condition of weak coupling, i.e., the coefficients
$C^1,\dots,C^N,K^1,\dots,K^N,L^1,\dots,L^N$ are sufficiently small, which makes $\mathcal{C},\mathcal{G}^X$ sufficiently small. The smallness of these coefficients de-emphasizes the role of the environment variables on the state dynamics as well as the cost functions of all DMs. This weak coupling condition along with the conditions for the stability of $\mathcal{A}$  guarantees the existence of $\Sigma\succeq0$ solving (\ref{eq:ssi2}).

The condition $\mbox{\rm cov}[X_0] = \Sigma$ is clearly a
strong one. A possible alternative setting
where steady-state initialization is not required would be the average cost setting.

\subsection{Case of identical DMs}
We formally introduce the notion of identical DMs as foolows.

\begin{defn}
\label{def:DM}
DM$^i$ and DM$^j$ are said to be identical if the following two conditions hold.
\begin{itemize}
\item[(i)] $(A_t^i,B_t^i,C_t^i,\beta^i,Q_t^i,R_t^i,K_t^i,L_t^i)=(A_t^j,B_t^j,C_t^j,\beta^j,Q_t^j,R_t^j,K_t^j,L_t^j)$, for all $t\in[0,T)$
\item[(ii)] $(x_0^i,(w_t^i)_{t\in[0,T)})$ and $(x_0^i,(w_t^j)_{t\in[0,T)})$ have identical distributions.
\end{itemize}
\end{defn}

When DM$^i$ and DM$^j$ are identical, $E^{\zeta^i} [J^i(s^i,\cdot)]=E^{\zeta^j} [J^j(s^j,\cdot)]$ for any $(s^i,\zeta^i)=(s^j,\zeta^j)$. 
We note that, due to the assumption $R_t^i\succ0$ for all $t\in[0,T)$,  a DM's optimal strategy with respect to any belief $\zeta\in\mathcal{P}(\mathbb{Y}^{[0,T)})$ is essentially unique, i.e., it can differ from that derived in Appendix only on a set of observation histories of probability zero.
Therefore, the optimal strategies of identical DMs with respect to any common belief $\zeta\in\mathcal{P}(\mathbb{Y}^{[0,T)})$ can differ only on a set of observation histories of probability zero. As a result, the strategies of identical DMs at an SEBEU (if exists) can be taken as identical.

We should also point out that, when DM$^i$ and DM$^j$ are identical, the estimates $E[x_k^i|Y_{k-1}]$ and $E[x_k^j|Y_{k-1}]$ generated by an SEBEU would be identical as well. Therefore, when DM$^i$ and DM$^j$ are identical, DMs can implement their SEBEU strategies by generating only one of $E[x_k^i|Y_{k-1}]$ or $E[x_k^j|Y_{k-1}]$.
For example, when all DMs are identical, the SEBEU strategies can be implemented by generating only $E[x_k^i|Y_{k-1}]$, for a generic DM$^i$, instead of generating 
the estimates $E[X_k|Y_{k-1}]$ of the entire system state.

\subsection{Case of countably infinite set of identical DMs}

We adopt the infinite-horizon setup of subsection~\ref{ss:lqgfp} with a countable set $(\textrm{DM}^i)_{i\in\mathbb{N}}$ of identical DMs.
The matrices $A^i$, $B^i$,  $Q^i$, $R^i$, $K^i$, $L^i$ are all assumed to be independent of the index $i$ (hence the superscript $i$ is dropped from these matrices), and the matrices $C^i$  are assumed to satisfy $C^i=0$, for all $i\in\mathbb{N}$, i.e., the individual state dynamics are not directly controlled through the environment variables.
Furthermore, the environment variable at each time $t\in\mathbb{N}_0$ is determined by the empirical average of the decisions as well as the states:
$$y_t = \limsup_{\bar{N} \to \infty} \frac{1}{\bar{N}} \sum_{j\in[1,\bar{N}]} (E^1 u_t^j + E^2 x_t^j) +\xi_t$$
where $E^1$, $E^2$ are constants of appropriate dimension.
In addition, we strengthen Assumption~\ref{as:lqudcinf} as follows.
\begin{assum}
\label{as:lqudcinfc}
$\mbox{}$

\begin{itemize}
\item[(i)] $(A,B)$ is stabilizable
\item[(ii)] $x_0^1,x_0^2,\dots$, $w_0^1,w_0^2,\dots$, $w_1^1,w_1^2,\dots$, $\xi_0,\xi_1,\dots$  are mutually independent  and have finite second-order moments
\item[(iii)] $(x_0^i)_{i\in\mathbb{N}}$, $(w_t^i)_{i\in\mathbb{N},t\in\mathbb{N}_0}$, $(\xi_t)_{t\in\mathbb{N}_0}$ each is iid.
\end{itemize}
\end{assum}

Consider an independent exogenous sequence $(z_t)_{t\in[0,\infty)}$ of environment variables which is iid with finite second-order moments and distribution $\zeta$.
As in the previous cases, suppose that each DM$^i$ instead aims to minimize $E^{\zeta} [J^i(s^i,\cdot)]$ over $s^i\in\mathbb{S}^i$. An optimal strategy for each DM$^i$ is obtained as
\begin{equation}
\label{eq:etawrz}
u_t^i=s_t^i(I_t^i,Z_{t-1}) = F x_t^i + G E[z_0] + H
\end{equation}
where $F$, $G:=\sum_{n\in\mathbb{N}_0} G_n$, $H$ are constants of appropriate dimensions,
which are pre-computable independently of the sequence $(z_t)_{t\in[0,\infty)}$ of environment variables; see Appendix~\ref{appx:lqinfinite}.

Consider now a joint strategy $s=(s^i)_{i\in\mathbb{N}} \in \mathbb{S}$ such that the sequence $(y_t)_{t\in[0,\infty)}$ of environment variables generated under $s$ is iid with distribution $\zeta_s$ ($E[|y_t|^2]=E[|y_0|^2]<\infty$, for all $t\in\mathbb{N}_0$). The joint strategy $s$ is an SEBEU if each $s^i$ minimizes $E^{\zeta_s} [J^i(s^i,\cdot)]$, e.g.,
\begin{align}
s_t^i(I_t^i,Y_{t-1}) =F x_t^i+ G E[y_0] + H  \label{eq:eqstrinf}
\end{align}
where
\begin{align}
y_t  = &  (E^1F+E^2) \limsup_{\bar{N}\rightarrow\infty}\frac{1}{\bar{N}}\sum_{j\in[1,\bar{N}]} x_t^j + E^1 G E[y_0] + E^1 H +\xi_t \label{eq:eqp2inf}\\
x_{t+1}^j  = & A x_t^j +B (F x_t^j + G E[y_0] + H ) + w_t^j, \  j\in\mathbb{N}. \label{eq:eqx2inf}
\end{align}
Therefore, an SEBEU can be constructed via (\ref{eq:eqstrinf}) provided an iid sequence \linebreak $(y_t)_{t\in[0,\infty)}$ of environment variables satisfying (\ref{eq:eqp2inf})-(\ref{eq:eqx2inf}) (as a random sequence)  can be obtained.
Due to Assumption~\ref{as:lqudcinfc}, $(x_t^i)_{i\in\mathbb{N}}$ generated by (\ref{eq:eqx2inf}) would be iid with finite second-order moments, for any $t\in\mathbb{N}_0$. As a result, we would have $\limsup_{\bar{N} \to \infty} \frac{1}{\bar{N}} \sum_{j\in[1,\bar{N}]} x_t^j=E[x_t^1]$, for any $t\in\mathbb{N}_0$, and
\begin{align*}
y_t  = &  (E^1 F + E^2) E[x_t^1] + E^1 G E[y_0] + E^1 H + \xi_t \\
E[x_{t+1}^1]  = & A E[x_t^1] +B (F E[x_t^1] + G E[y_0] + H )  + E[w_t^1].
\end{align*}
Thus, if the linear equations
\begin{align}
\hat{y}_0  = &  (E^1 F + E^2) \hat{x}_0 +  E^1 G \hat{y}_0 +  E^1 H + \hat{\xi}_0 \label{eq:hatp0} \\
\hat{x}_0  = & A \hat{x}_0 + B (F \hat{x}_0 + G \hat{y}_0 + H ) + \hat{w}_0^1. \label{eq:hatx0}
\end{align}
have a solution $(\hat{y}_0,\hat{x}_0)$ and $(E[y_0],E[x_0^1])=(\hat{y}_0,\hat{x}_0)$, then the environment variables generated by the strategies
\begin{equation}
\label{eq:etawrz3}
u_t^i =s_t^i(I_t^i,Y_{t-1})=F x_t^i + G \hat{y}_0 + H, \qquad i\in\mathbb{N}
\end{equation}
would be iid with finite second-order moments and satisfy (\ref{eq:eqp2inf})-(\ref{eq:eqx2inf}). This leads to the following result whose proof
follows from the preceding discussion.\\

\begin{thm}
\label{thm:exinfidm}
Consider the infinite-horizon linear quadratic problem of subsection~\ref{ss:lqgfp} with a countably infinite set of identical DMs under
Assumption~\ref{as:lqudcinfc}.
If $(E[y_0],E[x_0^1])=(\hat{y}_0,\hat{x}_0)$ where $(\hat{y}_0,\hat{x}_0)$ solves (\ref{eq:hatp0})-(\ref{eq:hatx0}),
there exists an SEBEU strategy $s\in\mathbb{S}$ defined by (\ref{eq:etawrz3}).\\
\end{thm}

In the case where the environment variables are determined by the average decisions of a countable set of identical DMs, no individual DM can change the sequence of environment variables by unilaterally deviating to an alternative strategy; therefore, an SEBEU is also a Nash equilibrium, which is in fact a mean-field equilibrium \cite{CainesMeanField2,lasry2007mean,jovanovic1988anonymous}.\\

%%%%%%%%%%%%%%%%%%%%%%%%%%%%%%%%%%%%%%%%%%%%%%%%%%%%%%%%%%%%%%%%%%%%%%%%%%%%%%%%

\section{An Example}

We present an example where $T=2$, all DMs are identical, and the environment variables are generated by 
\begin{align*}
y_t = &  \underbrace{\frac{1}{N} \sum_{i\in[1,N]} x_t^i}_{:= \bar{x}_t}  +\xi_t, \quad t\in[0,2).
\end{align*}
We also assume, in addition to Assumption~\ref{as:lqudc},  that all primitive random variables $$x_0^1,\dots,x_0^N,w_0^1,w_0^N,\dots,w_1^1,\dots,w_1^N,\xi_0,\xi_{1}$$ are mutually independent Gaussian random variables with zero-mean and
$\textrm{cov}[\xi_0]>0$.

In this case, the SEBEU established by Theorem~\ref{thm:ex} takes the form, for each $i\in[1,N]$,
\begin{align*}
s_0^i(x_0^i) = & F_0 x_0^i + G_{0,0} \hat{\bar{x}}_0 + G_{0,1} \hat{\bar{x}}_1\\
s_1^i(x_1^i,y_0)  = & F_1 x_1^i +G_{1,1} \hat{\bar{x}}_{1|0}(y_0)
\end{align*}
where $F_0,F_1,G_{0,0},G_{0,1},G_{1,1}$, are given in Appendix~\ref{appx:lqfinite} and 
\begin{align*}
\hat{\bar{x}}_0 & = E[\bar{x}_0]  \\ 
\hat{\bar{x}}_1 & = (A_0 + B_0 F_0) \hat{\bar{x}}_0 + B_0 (G_{0,0} \hat{\bar{x}}_0 + G_{0,1} \hat{\bar{x}}_1) + C_0 \hat{\bar{x}}_0\\
\hat{\bar{x}}_{1|0}(y_0) & = (A_0 + B_0 F_0) \hat{\bar{x}}_{0|0}(y_0)  + B_0 (G_{0,0} \hat{\bar{x}}_0 + G_{0,1}  \hat{\bar{x}}_1) + C_0y_0 \\
\hat{\bar{x}}_{0|0}(y_0) & = E[\bar{x}_0|y_0]= \textrm{cov}[\bar{x}_0]  (\textrm{cov}[\bar{x}_0]+\textrm{cov}[\xi_0])^{-1} y_0. 
\end{align*}
Since $E[\bar{x}_0]=0$ (by assumption) and 
\begin{equation}
\label{eq:IB0G0}
(I-B_0 G_{0,1}) \hat{\bar{x}}_1  = (A_0 + B_0 F_0 + B_0 G_{0,0} + C_0)E[\bar{x}_0] 
\end{equation}
an SEBEU exists corresponding to any $\hat{\bar{x}}_1$ in the null space of $I-B_0 G_{0,1}$. On the other hand, if $\bar{x}_0$ was not assumed to be zero-mean and (\ref{eq:IB0G0}) does not have a solution for $\hat{\bar{x}}_1$, an SEBEU does not exist.

Let us now consider the case where $I-B_0 G_{0,1}$ is non-singular (and $E[\bar{x}_0]=0$ as assumed). The SEBEU, which is essentially unique, is
\begin{align}
s_0^i(x_0^i) = & F_0 x_0^i \label{eq:ex3se0}\\
s_1^i(x_1^i,y_0)  = & F_1 x_1^i +G_1^N y_0 \label{eq:ex3se1}
\end{align}
where 
\begin{align*}
G_1^N  = & G_{1,1} \big((A_0 + B_0 F_0) \textrm{cov}[\bar{x}_0]  (\textrm{cov}[\bar{x}_0]+\textrm{cov}[\xi_0])^{-1}    + C_0\big).
\end{align*}
Due to $\textrm{cov}[\bar{x}_0]=\frac{1}{N}\textrm{cov}[x_0^1]$, we have $\lim_{N\rightarrow\infty} G_1^N=G_{1,1}C_0$; the state feedback gains $F_0,F_1$ are independent of the number of DMs.
Obtaining the set of Nash equilibria is not as straightforward; however, it is possible to show that, for any $\epsilon>0$, the SEBEU (\ref{eq:ex3se0})-(\ref{eq:ex3se1}) constitutes an $\epsilon-$Nash equilibrium if the number of DMs is sufficiently large.
To show this, we obtain DM$^i$'s optimal response $\tilde{s}^{i,N}\in\mathbb{S}^i$ to the strategies (\ref{eq:ex3se0})-(\ref{eq:ex3se1}) as follows.
\begin{align*}
\tilde{V}_1^i(I_1^i,y_0) := & \min_{u_1^i} \big\{ |x_1^i|_{Q_1}^2+|u_1^i|_{R_1}^2 +2  E[y_1^{\prime}| I_1^i,y_0] (K_1 u_1^i+L_1 x_1^i)  \\ & \qquad \ + \beta E \big[ |A_1 x_1^i + B_1 u_1^i + C_1 y_1 + w_1^i|_{Q_2}^2 | I_1^i,y_0 \big] \big \}.
\end{align*}
DM$^i$'s optimal response at stage $t=1$ is
$$\tilde{s}_1^{i,N}(I_1^i,y_0)=-(R_1 + \beta B_1^{\prime} Q_2 B_1)^{-1} (\beta B_1^{\prime} Q_2 A_1 x_1^i + (K_1^{\prime}+\beta B_1^{\prime} Q_2 C_1) E[y_1|I_1^i,y_0]). $$
We write
\begin{align*}
E[y_1|I_1^i,y_0] & = \frac{1}{N} x_1^i+(A_0+B_0 F_0)E[\bar{x}_0^{-i} | I_1^i,y_0]+\frac{N-1}{N} C_0y_0 
\end{align*}
where $\bar{x}_0^{-i}:=\frac{1}{N} \sum_{j\not=i}  x_0^j$.
From $y_0= \frac{1}{N}x_0^i + \bar{x}_0^{-i} + \xi_0$, we have
$$E[\bar{x}_0^{-i}|I_1^i,y_0]=\textrm{cov}[\bar{x}_0^{-i}] (\textrm{cov}[\bar{x}_0^{-i}]+\textrm{cov}[\xi_0])^{-1} \left(y_0-\frac{1}{N}x_0^i \right).$$
Therefore,
$$\tilde{s}_1^{i,N}(I_1^i,y_0)=\tilde{F}_1^N x_1^i +\tilde{G}_1^N y_0 +\tilde{F}_{1,0}^N x_0^i $$
where
\begin{align*}
\tilde{F}_1^N  := & F_1+\mathcal{O}(1/N) \\
\tilde{G}_1^N     := & G_{1,1} C_0 +\mathcal{O}(1/N) \\
\tilde{F}_{1,0}^N     := & \mathcal{O}(1/N).
\end{align*}
Let
\begin{align*}
\tilde{V}_0^i(x_0^i) := & \min_{u_0^i} \big\{|x_0^i|_{Q_0}^2+|u_0^i|_{R_0}^2   +2E[y_0^{\prime} | x_0^i] (K_0u_0^i+L_0 x_0^i)+  \beta E[\tilde{V}_1^i (I_1^i,y_0) | x_0^i)] \big\}.
\end{align*}
DM$^i$'s optimal response at stage $t=0$ can be derived as
$$\tilde{s}_0^{i,N}(x_0^i)=\tilde{F}_0^N x_0^i$$
where
\begin{align*}
\tilde{F}_0^N := & -  \big(R_0 + \beta B_0^{\prime} M_1 B_0  +\mathcal{O}(1/N) \big)^{-1} \big( \beta B_0^{\prime} M_1   A_0  +\mathcal{O}(1/N) \big) \\
M_1 =& Q_1 + \beta A_1^{\prime} \left(Q_2 - \beta Q_2 B_1 (R_1 + \beta B_1^{\prime} Q_2 B_1)^{-1} B_1^{\prime} Q_2 \right) A_1. 
\end{align*}

It is straightforward to see that $\lim_{N\rightarrow\infty}(\tilde{F}_0^N,\tilde{F}_1^N,\tilde{G}_1^N,\tilde{F}_{1,0}^N)= (F_0,F_1,G_{1,1}C_0,0)$.
This implies that, for any $\epsilon>0$, the SEBEU (\ref{eq:ex3se0})-(\ref{eq:ex3se1}) is an $\epsilon-$Nash equilibrium if the number of DMs is sufficiently large.

We should also remark that the SEBEU strategies (\ref{eq:ex3se0})-(\ref{eq:ex3se1}) are {\it always} linear mappings of each DM$^i$'s information variables whereas the same does not always hold for Nash equilibria for such linear quadratic Gaussian models \cite[Chapter 7]{basols99}.

%%%%%%%%%%%%%%%%%%%%%%%%%%%%%%%%%%%%%%%%%%%%%%%%%%%%%%%%%%%%%%%%%%%%%%%%%%%%%%%%
\section{Conclusion}
We extended some of the results in our earlier work on the notion of subjective equilibrium under beliefs of exogenous uncertainty (SEBEU) to the linear-quadratic case where the players interact through some environment variables. We established the existence of an SEBEU in pure strategies under the condition that a linear system of algebraic equations are solvable. Our results are constructive and the players strategies at an SEBEU are affine functions of their local state observations and their estimates of the system state. In the Gaussian case, the estimates of the system state can be generated by a Kalman filter. We pointed out that an SEBEU is also a Nash equilibrium in the case of countable number of identical players. Showing that an SEBEU is an approximate Nash equilibrium when there is a large but finite number of players remains as a future research topic.

\bibliographystyle{plain}
\bibliography{references}
%\nocite{*}

\appendix

\section{}

We drop the superscript $i$ and derive an optimal response for a generic DM with respect to independent exogenous environment variables $(z_t)_{t\in[0,T)}$ satisfying $\sup_{t\in[0,T)}E[|z_t|^2]<\infty$. DM has the knowledge of the distribution of $(z_t)_{t\in[0,T)}$ and has access to the realizations $Z_{t-1}=(z_k)_{k\in[0,t)}$ at time $t\in[0,T)$ prior to choosing $u_t$.
For any random vector $\zeta_t$, we use the notation $\hat{\zeta}_t=E[\zeta_t]$ and $\hat{\zeta}_{t|k}=E[\zeta_t|Z_k]$.
Assumption~\ref{as:lqudc} holds in the finite-horizon case, whereas Assumption~\ref{as:lqudcinf} holds in the infinite-horizon case.
By completing the squares, we write DM's undiscounted cost at time $t\in[0,T)$ as
$$c_t(x_t,u_t,z_t)=\tilde{c}_t(x_t,u_t,z_t) - |z_t|_{L_tQ_t^{-1}L_t^{\prime}+K_tR_t^{-1}K_t^{\prime}}^2$$
where
$$\tilde{c}_t(x_t,u_t,z_t):=|x_t + Q_t^{-1}L_t^{\prime} z_t|_{Q_t}^2+|u_t + R_t^{-1}K_t^{\prime}z_t|_{R_t}^2.$$
Since $E\big[\sum_{t\in[0,T)} \beta^t |z_t|_{K_tQ_t^{-1}K_t^{\prime}+L_tR_t^{-1}L_t^{\prime}}^2\big]$ is finite and independent of DM's strategy, we assume that DM's undiscounted cost function at time $t\in[0,T)$ is $\tilde{c}_t$ without loss of generality.

DM's optimal response in the finite-horizon case can be obtained using the main result in
\cite{duncan2012discrete} which studies the finite-horizon linear quadratic control problem with arbitrary correlated noise; however, for completeness, we provide a brief derivation of DM's optimal response below.

\subsection{Finite-horizon case}
\label{appx:lqfinite}

We define DM's cost-to-go functions as, for all $k\in [0,T]$ and $(I_{k},Z_{k-1})\in\mathbb{I}_k\times\mathbb{Y}^{[0,k)}$,
\begin{align*}
V_{k}(I_{k},Z_{k-1}) := & \min_{s_k,\dots,s_{T-1}} E \Bigg [\sum_{t\in[k,T)} \beta^{t-k} \tilde{c}_t(x_t,u_t,z_t)  +   \beta^{T-k} c_T(x_T) \Big| I_{k},Z_{k-1} \Bigg],
\quad k \in [0,T)
\end{align*}
with $V_{T}(I_{T},Z_{T-1}) :=   c_T(x_T)$, where $u_n=s_n(I_n,Z_{n-1})$, $s_n:(\mathbb{I}_n,\mathbb{Y}^{[0,n)}) \to \mathbb{U}$, for $n\in[k,T)$.
These cost-to-go functions satisfy, for $k\in[0,T)$,
\begin{align}
\label{eq:ctgX}
V_{k}(I_{k},Z_{k-1}) = & \min_{s_k}  E \big[ \tilde{c}_k(x_k,u_k,z_k)   + \beta    V_{k+1}(I_{k+1},Z_{k}) | I_{k},Z_{k-1}  \big].
\end{align}
As an induction hypothesis, assume that $V_{k+1}$ has the following form, for $k\in[0,T)$,
\begin{align}
&& V_{k+1}(I_{k+1},Z_k) = |x_{k+1}|_{M_{k+1}}^2  + 2x_{k+1}^{\prime} N_{k+1}(Z_k)  + O_{k+1}(Z_k) \label{eq:indhyp}
\end{align}
where $M_{k+1}\succeq0$ is an appropriate dimensional matrix, $N_{k+1}$, $O_{k+1}$ are appropriate dimensional functions of $Z_k$. $V_T$ satisfies this hypothesis with $M_T =  Q_T$, $N_T =  0$, $O_T = 0$.
The minimizing control in (\ref{eq:ctgX}) is obtained as
%\begin{align*}
%0 = 2 K_k^{\prime} \hat{z}_{k|k-1} + 2 R_k u_k + 2 \beta (B_k)^{\prime} M_{k+1} (A_k x_k + B_k u_k + C_k \hat{z}_{k|k-1}+\hat{w}_k) +  2 \beta (B_k)^{\prime} \hat{N}_{k+1|k-1}
%\end{align*}
\begin{align}
u_k = s_k(I_k,Z_{k-1}):=F_k x_k+G_{k,k} \hat{z}_{k|k-1} +\tilde{H}_k(Z_{k-1}) \label{eq:optstrX}
\end{align}
where
\begin{align}
F_k = & -\beta S_k^{-1}  B_k^{\prime} M_{k+1} A_k \label{eq:Fk} \\
G_{k,k} = & - S_k^{-1}  (K_k^{\prime}+ \beta B_k^{\prime} M_{k+1} C_k) \label{eq:Gkk}\\
\tilde{H}_k(Z_{k-1}) = & -\beta S_k^{-1} B_k^{\prime} (\hat{N}_{k+1|k-1} + M_{k+1} \hat{w}_k)\\
S_k = & R_k+ \beta B_k^{\prime} M_{k+1} B_k. \label{eq:Sk}
\end{align}
Substituting the minimizing control (\ref{eq:optstrX})-(\ref{eq:Sk}) in (\ref{eq:ctgX})-(\ref{eq:indhyp}) results in
\begin{align}
\nonumber V_{k}(I_{k},Z_{k-1}) = &  |x_k|_{Q_k}^2 + |F_k x_k+G_{k,k} \hat{z}_{k|k-1}  +\tilde{H}_k(Z_{k-1})|_{R_k}^2
\\ \nonumber & +2\hat{z}_{k|k-1}^{\prime}(K_k(F_k x_k+G_{k,k} \hat{z}_{k|k-1}  +\tilde{H}_k(Z_{k-1}))+L_kx_k)
\\ \nonumber &      + 2E[|z_k|_{K_k Q_k^{-1} K_k^{\prime}+L_kR_k^{-1}L_k^{\prime}}^2 | Z_{k-1} ]
\\ \nonumber & + \beta E  [ |(A_k+B_kF_k)x_k+B_k (G_{k,k} \hat{z}_{k|k-1} +\tilde{H}_k(Z_{k-1}))+C_kz_k+w_k|_{M_{k+1}}^2   | Z_{k-1} ]
\\ \nonumber & + \beta E  [ 2((A_k+B_kF_k)x_k+B_k (G_{k,k} \hat{z}_{k|k-1} +\tilde{H}_k(Z_{k-1}))+C_kz_k+w_k)^{\prime} N_{k+1}(Z_k)  | Z_{k-1} ]
\\ & + \beta E  [  O_{k+1}(Z_k) | Z_{k-1} ]. \label{eq:Vk}
\end{align}
Grouping the constant, linear, and quadratic terms in $x_k$ for fixed $Z_{k-1}$, we rewrite $V_{k}(I_{k},Z_{k-1})$ as
\begin{align*}
V_{k}(I_{k},Z_{k-1}) = &  |x_{k}|_{M_{k}}^2 + 2x_k^{\prime} N_{k}(Z_{k-1})   + O_{k}(Z_{k-1})
\end{align*}
where
\begin{align}
\label{eq:ricL}
M_k = & Q_k + \beta A_k^{\prime} (M_{k+1} - \beta M_{k+1} B_k S_k^{-1}  B_k^{\prime} M_{k+1}) A_k \\
\label{eq:ricM}
N_k(Z_{k-1}) = &  \beta (A_k+B_kF_k)^{\prime} \hat{N}_{k+1|k-1} + \hat{v}_{k|k-1} \\
v_k = & (F_k^{\prime}K_k^{\prime} +L_k^{\prime}) z_k   + \beta (A_k+B_kF_k)^{\prime}  M_{k+1} (C_k z_k + w_k)
\\
O_k (Z_{k-1}) = & \beta \hat{O}_{k+1|k-1} + \hat{\lambda}_{k|k-1}\\
\lambda_k  = &
\nonumber |G_{k,k} \hat{z}_{k|k-1}  +\tilde{H}_k(Z_{k-1})|_{R_k}^2+
2z_k^{\prime}K_k(G_{k,k} \hat{z}_{k|k-1}  +\tilde{H}_k(Z_{k-1})) \\ & + |z_k|_{K_k Q_k^{-1} K_k^{\prime}+L_kR_k^{-1}L_k^{\prime}}^2
+ \nonumber \beta |B_k (G_{k,k} \hat{z}_{k|k-1} +\tilde{H}_k(Z_{k-1}))+C_kz_k+w_k|_{M_{k+1}}^2
\\ & + 2 \beta  (B_k (G_{k,k} \hat{z}_{k|k-1} +\tilde{H}_k(Z_{k-1}))+C_kz_k+w_k)^{\prime} N_{k+1}(Z_k) .
\label{eq:ricO}
\end{align}
The recursion (\ref{eq:ricL}) is the well-known Riccati difference equation for discrete time linear quadratic problems (without the environment variables), $N_k(Z_{k-1})$ is obtained from the right hand side of (\ref{eq:Vk}) by collecting the linear terms in $x_k$ for fixed $Z_{k-1}$ and using $R_k F_k+\beta B_k^{\prime} M_{k+1} (A_k+B_kF_k)=0$, and $O_k(Z_{k-1})$ equals the right hand side of (\ref{eq:Vk}) for $x_k=0$.
Since $V_{k}(I_{k},Z_{k-1})$ satisfies the induction hypothesis, an optimal strategy for DM is given by (\ref{eq:optstrX}). This optimal strategy is rewritten as, for $k\in[0,T)$,
\begin{align*}
s_k(I_k,Z_{k-1}) = & F_k x_k +  \sum_{t\in[k,T)} G_{k,t} \hat{z}_{t|k-1} + H_k
\end{align*}
where
\begin{align*}
G_{k,k} = & - S_k^{-1}  (K_k^{\prime}+ \beta B_k^{\prime} M_{k+1} C_k) \\
G_{k,t} = & -\beta S_k^{-1} B_k^{\prime} (\Phi_{k,t} M_{t+1} C_t +\Phi_{k,t-1} (F_t^{\prime}K_t^{\prime}+L_t^{\prime})), \quad t \in [k+1,T)     \\
H_k = & - \beta S_k^{-1}  B_k^{\prime} \sum_{t\in[k,T)} \Phi_{k,t} M_{t+1} \hat{w}_t \\
\Phi_{k,k} =& I, \quad \Phi_{k,t} =  \prod_{n\in[k+1,t]} \beta (A_{n} + B_{n} F_{n})^{\prime}, \ \ t\in[k+1,T).
\end{align*}

\subsection{Infinite-horizon case}
\label{appx:lqinfinite}

Recall that $A_t,B_t,C_t, K_t, L_t, Q_t, R_t$ are independent of $t$ (hence the subscript $t$ is dropped) and $\beta\in(0,1)$.

Let $J_k^*(x_k,Z_{k-1})$ denote DM's optimal cost for the infinite-horizon problem starting from $(x_k,Z_{k-1})$ at time $k\in\mathbb{N}_0$.
For any $\tilde{T}\in\mathbb{N}$, let $J_{k,\tilde{T}}^*(x_k,Z_{k-1})$ denote DM's optimal cost for the finite-horizon version of the problem starting from $(x_k,Z_{k-1})$ at time $k\in\mathbb{N}_0$ and running over the horizon $[k,\tilde{T})$ (with no terminal cost). We have, for $k\in\mathbb{N}_0$, $$J_{k,k+1}^*(x_k,Z_{k-1}) \leq J_{k,k+2}^*(x_k,Z_{k-1}) \leq \cdots \leq J_k^*(x_k,Z_{k-1})<\infty$$
where $J_k^*(x_k,Z_{k-1})<\infty$ is due to the stabilizability of $(A,B)$ and $\sup_{t\in\mathbb{N}_0}E[|z_t|^2]<\infty$. Hence, for any $k\in\mathbb{N}_0$ and $(x_k,Z_{k-1})$, $\lim_{\tilde{T}\rightarrow\infty} J_{k,\tilde{T}}^*(x_k,Z_{k-1})=J_{k,\infty}^*(x_k,Z_{k-1})$ for some $J_{k,\infty}^*(x_k,Z_{k-1}) \leq J_k^*(x_k,Z_{k-1})$. Note that
$$
J_{k,\tilde{T}}^*(x_k,Z_{k-1})= |x_k|_{M_k}^2 + 2x_k^{\prime} N_k(Z_{k-1})   + O_k(Z_{k-1})
$$
where $M_k, N_k(Z_{k-1}), O_k(Z_{k-1})$, $k\in\mathbb{N}_0$, are generated by the backward recursions (\ref{eq:ricL})-(\ref{eq:ricO}) with the boundary conditions $M_{\tilde{T}}=N_{\tilde{T}}=O_{\tilde{T}}=0$.
This implies that, for any fixed $k\in\mathbb{N}_0$ and $Z_{k-1}$, $M_k$, $N_k(Z_{k-1})$, $O_k(Z_{k-1})$
converge to some finite values as $\tilde{T}\rightarrow\infty$.

It is well known that, as $\tilde{T}\rightarrow\infty$, $M_k$ for any fixed $k\in\mathbb{N}_0$ tends to the unique solution $M\succeq0$ of the algebraic Riccati equation 
\begin{align}
\label{eq:riccL}
M = & Q + \beta A^{\prime} (M - \beta M B S^{-1}  B^{\prime} M) A
\end{align}
within the set of positive semi-definite matrices since $(A,B)$ is stabilizable and $(A,Q^{1/2})$ is detectable.
Therefore, as $\tilde{T}\rightarrow\infty$,  $F_k$, $G_{k,k}$, $S_k$ for any fixed $k\in\mathbb{N}_0$ introduced in (\ref{eq:Fk})-(\ref{eq:Sk})
converge to
\begin{align*}
F := & -\beta S^{-1}  B^{\prime} M A  \\
G := & - S^{-1}  (K^{\prime}+ \beta B^{\prime} M C) \\
S := & R+ \beta B^{\prime} M B
\end{align*}
respectively, where $F$ renders the matrix $\sqrt{\beta} (A+BF)$  asymptotically stable. It follows from (\ref{eq:ricL})-(\ref{eq:ricO}) that, as $\tilde{T}\rightarrow\infty$, $N_k(Z_{k-1})$, $\tilde{H}_k(Z_{k-1})$, $O_k(Z_{k-1})$  for any fixed $k\in\mathbb{N}_0$ and $Z_{k-1}$ tend to
\begin{align*}
N_{k,\infty}(Z_{k-1}) := & \sum_{n\in\mathbb{N}_0 } (\beta (A+BF)^{\prime})^n \hat{\bar{v}}_{k+n|k-1} \\
\tilde{H}_{k,\infty}(Z_{k-1}) := & -\beta S^{-1} B^{\prime} \bigg(\sum_{n\in\mathbb{N}_0 } (\beta (A+BF)^{\prime})^n \hat{v}_{k+n+1|k-1} + M \hat{w}_0\bigg) \\
O_{k,\infty}(Z_{k-1}) = & \sum_{n\in\mathbb{N}_0 } \beta^n \hat{\bar{\lambda}}_{k+n|k-1}.
\end{align*}
respectively, where
\begin{align*}
\bar{v}_k := &  (F^{\prime}K^{\prime} +L) z_k   + \beta (A+BF)^{\prime}  M (C z_k + w_k)\\
\bar{\lambda}_k  := &
|G \hat{z}_{k|k-1}  +\tilde{H}_{k,\infty}(Z_{k-1})|_{R}^2+
2z_k^{\prime}K(G \hat{z}_{k|k-1}  +\tilde{H}_{k,\infty}(Z_{k-1})) \\ & + |z_k|_{K Q^{-1} K^{\prime}+LR^{-1}L^{\prime}}^2
+ \beta |B (G \hat{z}_{k|k-1} +\tilde{H}_{k,\infty}(Z_{k-1}))+Cz_k+w_k|_{M}^2
\\ & + 2 \beta  (B (G \hat{z}_{k|k-1} +\tilde{H}_{k,\infty}(Z_{k-1}))+Cz_k+w_k)^{\prime} N_{k+1,\infty}(Z_k).
\end{align*}

Straightforward computation shows that $(J_{k,\infty})_{k\in\mathbb{N}_0}$ satisfy the Bellman equations, i.e.,
\begin{align}
\nonumber J_{k,\infty}^*(x_k,Z_{k-1}) & = |x_k|_M^2 + 2x_k^{\prime} N_{k,\infty}(Z_{k-1})   + O_{k,\infty}(Z_{k-1}) \\
& = \min_{u_k\in\mathbb{U}} E[|x_k + Q^{-1}L^{\prime} z_k|_{Q}^2+|u_k + R^{-1}K^{\prime}z_k|_{R}^2 + \beta J_{k+1,\infty}^*(x_{k+1},Z_k) | x_k,Z_{k-1}] \label{eq:bellman}
\end{align}
with the minimizing control
\begin{align}
 u_k^*  & =   s_k^*(x_k,Z_{k-1}):=F x_k +G\hat{z}_{k|k-1} +\tilde{H}(Z_{k-1}). \label{eq:bellman2}
\end{align}
Now, let $J_{k,s^*}(x_k,Z_{k-1})$ denote the cost achieved by the strategy $s^*=(s_0^*,s_1^*,\dots)$ starting from $(x_k,Z_{k-1})$ at time $k\in\mathbb{N}_0$.
Let $(x_{n+1}^*,u_n^*)_{n\in[k,\tilde{T})}$ be generated by $(s_n^*,z_n)_{n\in[k,\tilde{T})}$ from $(x_k^*=x_k,Z_{k-1})$.
Since $(J_{k,\infty}^*)_{k\in\mathbb{N}_0}$ is non-negative valued, we write
\begin{align*}
J_{k,s^*}(x_k,Z_{k-1})
\leq & \lim_{\tilde{T} \rightarrow \infty} E \Bigg [\sum_{n\in[k,\tilde{T})} \beta^{n-k} \tilde{c}_n(x_n^*,u_n^*,z_n)+\beta^{\tilde{T}-k} J_{\tilde{T},\infty}^*(x_{\tilde{T}}^*,Z_{\tilde{T}-1}) \bigg| x_k,Z_{k-1} \Bigg]
\\
= & \lim_{\tilde{T} \rightarrow \infty} E \Bigg [\sum_{n\in[k,\tilde{T}-1)} \beta^{n-k} \tilde{c}_n(x_n^*,u_n^*,z_n)+\beta^{\tilde{T}-k-1} J_{\tilde{T}-1,\infty}^*(x_{\tilde{T}-1}^*,Z_{\tilde{T}-2})  \bigg| x_k,Z_{k-1} \Bigg]7
\\ \vdots & \\
= &  J_{k,\infty}^*(x_k,Z_{k-1})
\end{align*}
where the equalities are obtained by repeated application of (\ref{eq:bellman})-(\ref{eq:bellman2}). Therefore, the strategy $s^*$ is optimal.

Finally, consider any optimal strategy $\check{s}=(\check{s}_0,\check{s}_1,\dots)$. We must have, for $k\in\mathbb{N}_0$,
\begin{align*}
J_k^*(x_k,Z_{k-1}) & = E[\tilde{c}_k(x_k,\check{u}_k,z_k) + \beta J_{k+1}^*(Ax_k+B\check{u}_k+Cz_k+w_k,Z_k) | x_k,Z_{k-1}] \\
 & = \min_{u_k} E[\tilde{c}_k(x_k,u_k,z_k) + \beta J_{k+1}^*(Ax_k+Bu_k+Cz_k+w_k,Z_k) | x_k,Z_{k-1}]
\end{align*}
where $\check{u}_k=\check{s}_k(x_k,Z_{k-1})$
Since the minimum above is uniquely achieved by $s_k^*(x_k,Z_{k-1})$, we have $\check{s}=s^*$.
The optimal strategy is rewritten as
\begin{align*}
u_k^* = s_k^*(x_k,Z_{k-1}) = F x_k  +\sum_{n\in\mathbb{N}_0} G_n \hat{z}_{k+n|k-1} + H
\end{align*}
where
\begin{align*}
G_0 = & - S^{-1} (K^{\prime}+\beta B^{\prime} M C) \\
G_n = & -\beta S^{-1} B^{\prime} (\beta (A+BF)^{\prime})^n    ((F^{\prime}K^{\prime} +L^{\prime}) + \beta (A+BF)^{\prime}  M C  ), \quad n \in\mathbb{N}    \\
H = & -\beta S^{-1}  B^{\prime} (I - \beta (A + B F)^{\prime})^{-1}  M \hat{w}_0 .
\end{align*}

\end{document}